\documentclass[reqno]{amsart}

\usepackage{amssymb}

\usepackage[dvips,colorlinks=true,linkcolor=blue]{hyperref}

\newcounter{numcount}
\newcommand{\labelnummer}{\mbox{\normalfont (\roman{numcount})}}%

\makeatletter

   {\let\curlabelspeicher\@currentlabel%
     \begin{list}{\labelnummer}%
       {\usecounter{numcount}\leftmargin0pt%
         \topsep0.5ex\partopsep2ex\parsep0pt\itemsep0ex\@plus1\p@%
         \labelwidth2.5em\itemindent3.5em\labelsep1em%
       }%
     \let\saveitem\item%
     \def\item{\saveitem%
  \def\@currentlabel{\curlabelspeicher\hskip0.5pt\labelnummer}}%
     \let\savelabel\label%
     \def\label##1{\savelabel{##1}%
       \@bsphack%
         \ifmmode\else%
           \protected@write\@auxout{}%
           {\string\newlabel{##1item}{{\labelnummer}{\thepage}}}%
         \fi%
       \@esphack%
     }%
   }{\end{list}}%

\newcommand{\esp}{\mathbb{E}}

\newcommand{\vers}{\operatornamewithlimits{\to}}
\theoremstyle{plain}
\newtheorem{Th}{Theorem}[section]
\newtheorem{Le}{Lemma}[section]
\newtheorem{Pro}{Proposition}[section]

\theoremstyle{definition}
\newtheorem{Rem}{Remark}[section]

\newcommand\R{\mathbb R}
\newcommand\N{\mathbb N}

\newcommand\Z{\mathbb Z}

\renewcommand\P{\mathbb P}

\newcommand\e{\mathrm{e}}

\newcommand\beq{\begin{equation}}
\newcommand\eeq{\end{equation}}

\begin{document}

\title[Spectral statistics in the localized regime]{Spectral
  statistics for the discrete Anderson model in the localized
  regime} 

\author{Fran\c cois Germinet} \address{Universit{\'e} de Cergy-Pontoise,
  CNRS UMR 8088, IUF, D{\'e}partement de math{\'e}matiques, F-95000
  Cergy-Pontoise, France}
\email{\href{mailto:francois.germinet@u-cergy.fr}{francois.germinet@u-cergy.fr}}
\author{Fr{\'e}d{\'e}ric Klopp} \address{LAGA, Institut Galil{\'e}e, CNRS UMR
  7539, Universit{\'e} de Paris-Nord, Avenue J.-B.  Cl{\'e}ment, F-93430
  Villetaneuse, France\ and \ Institut Universitaire de France}
\email{\href{mailto:klopp@math.univ-paris13.fr}{klopp@math.univ-paris13.fr}}


\thanks{This text is a contribution to the proceedings of the
  conference ``Spectra of Random Operators and Related Topics'' held
  at Kyoto University, 02-04/12/09 organized by N. Minami and
  N. Ueki. \\ Both authors are supported by the grant
  ANR-08-BLAN-0261-01. \\ Part of this work was done while the authors
  were guests at Centre Interfacultaire Bernoulli (EPFL). Both authors
  acknoledge support from the FNS}

%

\maketitle

\begin{abstract} 
  We report on recent results on the spectral statistics of the
  discrete Anderson model in the localized phase obtained
  in~\cite{Ge-Kl:10}. In particular, we describe the
  \begin{itemize}
  \item locally uniform Poisson behavior of the rescaled eigenvalues,
  \item independence of the Poisson processes obtained as such limits
    at distinct energies,
  \item locally uniform Poisson behavior of the joint distributions of
    the rescaled energies and rescaled localization centers in a large
    range of scales.
  \item the distribution of the rescaled level spacings, locally and
    globally in energy,
  \item the distribution of the rescaled localization centers
    spacings.
  \end{itemize}
  Our results show, in particular, that, for the discrete Anderson
  Hamiltonian with smoothly distributed random potential at
  sufficiently large coupling, the limit of the level spacing
  distribution is that of i.i.d. random variables distributed
  according to the density of states of the random Hamiltonian.
\end{abstract}

\section{Introduction}
\label{sec:introduction}
On $\ell^2(\Z^d)$, consider the random Anderson model
\begin{equation*}
  H_{\omega}=-\Delta+V_\omega
\end{equation*}
where $-\Delta$ is the discrete Laplace operator
\begin{equation*}
  (-\Delta u)_n=\sum_{|m-n|=1}u_m\quad\text{ for }
  u=(u_n)_{n\in\Z^d}\in\ell^2(\Z^d)
\end{equation*}
and $V_\omega$ is the random potential
\begin{equation*}
  (V_\omega u)_n=\omega_n u_n\quad\text{ for }u=(u_n)_{n\in\Z^d}
  \in\ell^2(\Z^d).
\end{equation*}
We assume that the random variables $(\omega_n)_{n\in\Z^d}$ are
independent identically distributed and that their distribution admits
a compactly supported bounded density, say $g$.\\
It is then well known (see e.g.~\cite{MR2509110}) that
\begin{itemize}
\item there exists $\Sigma:=[S_-,S_+]=[-2d,2d]+$supp$\,g\subset\R$
  such that, for almost every $\omega=(\omega_n)_{n\in\Z^d}$, the
  spectrum of $H_{\omega}$ is equal to $\Sigma$;
\item for some $S_-<s_-\leq s_+<S_+$, the intervals $I_-=[S_-,s_-)$
  and $I_+=(s_+,S_+]$ are contained in the region of complete
  localization for $H_\omega$, in particular, $I_-\cup I_+$ contains
  only pure point spectrum associated to exponentially decaying
  eigenfunctions; for the precise meaning of the region of complete
  localization, we refer to~\cite{MR2002h:82051,MR2509110,MR2203782};
  if the disorder is sufficiently large or if the dimension $d=1$
  then, one can pick $I_+\cup I_-=\Sigma$; define $I=I_+\cup I_-$;
\item there exists a bounded density of states, say
  $E\mapsto\nu(E)$, such that, for any continuous function
  $\varphi:\ \R\to\R$, one has
  \begin{equation}
    \label{eq:10}
    \int_\R\varphi(E)\nu(E)dE=
    \mathbb{E}(\langle\delta_0,\varphi(H_\omega)\delta_0\rangle).
  \end{equation}
  Here, and in the sequel, $\mathbb{E}(\cdot)$ denotes the expectation
  with respect to the random parameters.\\
  Let $N$ be the integrated density of states of $H_\omega$ i.e. $N$
  is the distribution function of the measure $\nu(E)dE$. The function
  $\nu$ is only defined $E$ almost everywhere. In the sequel, unless
  we explicitly say otherwise, when we speak of $\nu(E)$ for some $E$,
  we mean that the non decreasing function $N$ is differentiable at
  $E$ and that $\nu(E)$ is its derivative at $E$.
\end{itemize}
We now describe the local level and localization center statistics,
the level spacing statistics and the localization center spacings
statistics in $I$.
\section{The local level statistics}
\label{sec:level-statistics}
For $L\in\N$, let $\Lambda=\Lambda_L=[-L,L]^d\cap\Z^d\subset\Z^d$ be a
large box and $H_{\omega,\Lambda}$ be the operator $H_\omega$
restricted to $\Lambda$ with periodic boundary conditions. Let $N$ be
the volume of $\Lambda$ i.e. $N=(2L+1)^d$.
\par $H_\omega(\Lambda)$ is an $N\times N$ real symmetric matrix. Let us
denote its eigenvalues ordered increasingly and repeated according to
multiplicity by $E_1(\omega,\Lambda)\leq E_2(\omega,\Lambda)\leq
\cdots\leq E_N(\omega,\Lambda)$.
\par Let $E_0$ be an energy in $I$ such that $\nu(E_0)>0$.  The local
level statistics near $E_0$ is the point process defined by
\begin{equation}
  \label{eq:13}
  \Xi(\xi,E_0,\omega,\Lambda) = 
  \sum_{j=1}^N \delta_{\xi_j(E_0,\omega,\Lambda)}(\xi)
\end{equation}
where
\begin{equation}
  \label{eq:11}
  \xi_j(E_0,\omega,\Lambda)=|\Lambda|\,\nu(E_0)\,(E_j(\omega,\Lambda)-
  E_0),\quad 1\leq j\leq N.
\end{equation}
The main result of~\cite{MR97d:82046} reads
\begin{Th}[\cite{MR97d:82046}]
  \label{thr:3}
  Let $E_0$ be an energy in $I$ such that $\nu(E_0)>0$. When
  $|\Lambda|\to+\infty$, the point process $\Xi(E_0,\omega,\Lambda)$
  converges weakly to a Poisson process on $\R$ with intensity the
  Lebesgue measure.
\end{Th}
\subsection{Uniform Poisson convergence}
\label{sec:strong-poiss-conv}
In~\cite{Ge-Kl:10}, we obtain a uniform version of Theorem~\ref{thr:3}
i.e. a version that holds uniformly over an energy interval of size
asymptotically infinite compared to $|\Lambda|^{-1}$.\\
Fix $1>\beta>d/(d+2)$.  Let $I_\Lambda(E_0,\beta)$ be the interval
centered at $E_0$ of length $2|\Lambda|^{-\beta}$. Let the number of
eigenvalues of $H_\omega(\Lambda)$ inside $I_\Lambda(E_0,\beta)$ be
equal to $N_\Lambda(\omega,E_0)$.  For $1\leq j\leq
N_\Lambda(\omega,E_0)-1$, define the renormalized eigenvalues
$\xi_j(\omega,\Lambda)$ by~\eqref{eq:11} for $E_j\in
I_\Lambda(E_0,\beta)$.  Hence, for all $1\leq j\leq
N_\Lambda(\omega,E_0)-1$, one has
$\xi_j(\omega,\Lambda)\in|\Lambda|^{1-\beta}\cdot[-1,1]$.\\
We then prove
\begin{Th}[\cite{Ge-Kl:10}]
  \label{thr:2}
  Let $E_0$ be an energy in $I$ such that $\nu(E_0)>0$.\\
  Then, there exists $\delta>0$, such that, for any sequences of
  intervals $I_1=I^\Lambda_1,\dots,I_p=I^\Lambda_p$ in
  $|\Lambda|^{1-\beta}\cdot[-1,1]$ such that
  \begin{equation}
    \label{eq:32}
    \inf_{j\not= k}\text{dist}(I_j, I_k)\geq e^{-|\Lambda|^\delta},
  \end{equation}
  one has, for any sequences of integers
  $k_1=k^\Lambda_1,\cdots,k_p=k^\Lambda_p\in\N^p$,
  \begin{equation*}
    \lim_{|\Lambda|\to+\infty}
    \left|\P\left(\left\{\omega;\
        \begin{split}
          &\#\{j;\ \xi_j(\omega,\Lambda)\in
          I_1\}=k_1\\&\vdots\hskip3cm\vdots\\
          &\#\{j;\ \xi_j(\omega,\Lambda)\in I_p\}=k_p
        \end{split}
      \right\}\right)-e^{-|I_1|}\frac{|I_1|^{k_1}}{k_1!}\cdots
    e^{-|I_p|}\frac{|I_p|^{k_p}}{k_p!}\right|=0.
  \end{equation*}
\end{Th}
\noindent Note that, in Theorem~\ref{thr:2}, we don't require the
limits
\begin{equation*}
  \begin{split}
    &\lim_{|\Lambda|\to+\infty}e^{-|I_1|}\frac{|I_1|^{k_1}}{k_1!}=
    \lim_{|\Lambda|\to+\infty}e^{-|I^\Lambda_1|}\frac{|I^\Lambda_1|^{k^\Lambda_1}}{k^\Lambda_1!},\    \dots,\\
    &\lim_{|\Lambda|\to+\infty}e^{-|I_p|}\frac{|I_p|^{k_p}}{k_p!}=
    \lim_{|\Lambda|\to+\infty}e^{-|I^\Lambda_p|}\frac{|I^\Lambda_p|^{k^\Lambda_p}}{k^\Lambda_p!}
  \end{split}
\end{equation*}
to exist. \\
Clearly, Theorem~\ref{thr:3} is a consequence of the stronger
Theorem~\ref{thr:2}. The main improvement over the statements found
in~\cite{MR97d:82046} is that the interval over which the Poisson
statistics holds uniformly is much larger. We also note that
Theorem~\ref{thr:2} gives the asymptotics of the level spacing
distribution over intervals $I_\Lambda$ of size $|\Lambda|^{-d/(d+2)}$
(see section~\ref{sec:level-statistics} and, in particular,
Theorem~\ref{thr:1}). It also gives the asymptotic independence of the
local Poisson processes defined at energies $E_\Lambda$ and
$E'_\Lambda$ such that
\begin{equation*}
 |E_\Lambda-E_0|+|E'_\Lambda-E_0|\leq
  |\Lambda|^{-\beta}\quad\text{ and }\quad
\displaystyle |\Lambda|\cdot|E_\Lambda-E'_\Lambda|
  \vers_{\Lambda\to\Z^d}+\infty  
\end{equation*}
We refer to the next section for more general results on this
asymptotic independence.\\
It is natural to wonder what is the largest size of interval in which
a result like Theorem~\ref{thr:2}. We do not know the answer to that
question.
\subsection{Asymptotic independence of the local processes}
\label{sec:asympt-indep-local}
Once Theorem~\ref{thr:3} is known, it is natural to wonder how the
point processes obtained at two distinct energies relate to each
other. We prove the following
\begin{Th}[\cite{Ge-Kl:10,Kl:10}]
  \label{thr:4}
  Assume that the dimension $d=1$. Pick $E_0\in I$ and $E'_0\in I$
  such that $E_0\not=E'_0$, $\nu(E_0)>0$ and $\nu(E'_0)>0$.\\
  When $|\Lambda|\to+\infty$, the point processes
  $\Xi(E_0,\omega,\Lambda)$ and $\Xi(E'_0,\omega,\Lambda)$, defined
  in~\eqref{eq:13}, converge weakly respectively to two independent
  Poisson processes on $\R$ with intensity the Lebesgue measure. That
  is, for $U_+\subset\R$ and $U_-\subset\R$ compact intervals and
  $\{k_+,k_-\}\in\N\times\N$, one has
  \begin{equation*}
    \mathbb{P}\left(\left\{\omega;\
        \begin{cases}
          &\#\{j;\xi_j(E_0,\omega,\Lambda)\in U_+\}=k_+\\
          &\#\{j;\xi_j(E'_0,\omega,\Lambda)\in U_-\}=k_-
        \end{cases}
      \right\}
    \right)\vers_{\Lambda\to\Z^d}e^{-|U_+|}\frac{|U_+|^{k_+}}{k_+!}\cdot
    e^{-|U_-|}\frac{|U_-|^{k_-}}{k_-!}.
  \end{equation*}
\end{Th}
\noindent So we see that, in the localized regime, in dimension 1, at
distinct energies, the local eigenvalues behave independently from
each other. Theorem~\ref{thr:4} is a consequence of a decorrelation
estimate for distinct energies that is proved in~\cite{Kl:10}. It is
natural to expect that this decorrelation estimate stays true and,
hence, that Theorem~\ref{thr:4} stays true for arbitrary
dimensions. Nevertheless, we are only able to prove
\begin{Th}[\cite{Ge-Kl:10,Kl:10}]
  \label{thr:18}
  Pick $E_0\in I$ and $E'_0\in I$
  such that $|E_0-E'_0|>2d$, $\nu(E_0)>0$ and $\nu(E'_0)>0$.\\
  When $|\Lambda|\to+\infty$, the point processes
  $\Xi(E_0,\omega,\Lambda)$ and $\Xi(E'_0,\omega,\Lambda)$, defined
  in~\eqref{eq:13}, converge weakly respectively to two independent
  Poisson processes on $\R$ with intensity the Lebesgue measure.
\end{Th}
\noindent Theorems~\ref{thr:4} and~\ref{thr:18} naturally lead to
wonder how far the energies $E_0$ and $E'_0$ need to be from each
other with respect to the scaling used to renormalize the eigenvalues
for the asymptotic independence to still hold.\\
We prove
\begin{Th}[\cite{Ge-Kl:10}]
  \label{thr:8}
  Pick $E_0\in I$ such that $\nu(E_0)>0$. Assume moreover that the
  density of states $\nu$ is continuous at $E_0$.\\
  Consider two sequences of energies, say $(E_\Lambda)_\Lambda$ and
  $(E'_\Lambda)_\Lambda$ such that
  \begin{enumerate}
  \item one has $\displaystyle E_\Lambda\vers_{\Lambda\to\Z^d}E_0$ and
    $\displaystyle E'_\Lambda\vers_{\Lambda\to\Z^d}E_0$,
  \item one has $\displaystyle |\Lambda|\cdot|E_\Lambda-E'_\Lambda|
    \vers_{\Lambda\to\Z^d}+\infty$.
  \end{enumerate}
  Then, the point processes $\Xi(E_\Lambda,\omega,\Lambda)$ and
  $\Xi(E'_\Lambda,\omega,\Lambda)$, defined in~\eqref{eq:13}, converge
  weakly respectively to two independent Poisson processes on $\R$
  with intensity the Lebesgue measure.
\end{Th}
\noindent A crucial tool in proving Theorem~\ref{thr:8} are the
generalized Minami estimates proved in~\cite{MR2505733} that can also
be interpreted as local decorrelation estimates. Theorem~\ref{thr:8}
shows that, in the localized regime, eigenvalues that are sufficiently
far away from each other but still close, i.e. that are separated by a
distance that is asymptotically infinite when compared to the mean
spacing between the eigenlevels, behave as independent random
variables. There are no interactions except at very short
distances.\\
Assumption (2) can clearly not be omitted in Theorem~\ref{thr:8}; it
suffices to consider e.g. $E_\Lambda=E'_\Lambda+ a|\Lambda|^{-1}$ to
see that the two limit random processes are obtained as a shift from
one another. \\
To complete this section, we note again that, when
$|E'_\Lambda-E_\Lambda|=o(|\Lambda|^{-d/(d+2)})$, Theorem~\ref{thr:8}
is a consequence of Theorem~\ref{thr:2}.
\section{Localization center statistics}
\label{sec:local-cent-stat}
Recall that $E_1(\omega,\Lambda)\leq E_2(\omega,\Lambda)\leq
\cdots\leq E_N(\omega,\Lambda)$ denote the eigenvalues of
$H_{\omega,\Lambda}$ ordered increasingly and repeated according to
multiplicity.
\par To $E_j(\omega,\Lambda)$, we associate a normalized eigenvector
of $H_{\omega,\Lambda}$, say $\varphi_j(\omega,\Lambda)$. The
components of the vector $\varphi_j(\omega,\Lambda)$ are denoted by
$(\varphi_j(\omega,\Lambda; \gamma))_{\gamma\in\Lambda}$.\\
For $\varphi\in\ell^2(\Lambda)$, define the set of localization
centers for $\varphi$ as
\begin{equation*}
  C(\varphi)=\{\gamma\in\Lambda;\
  \varphi(\gamma)=\max_{\gamma'\in\Lambda}|\varphi(\gamma')|\}.
\end{equation*}
One has
\begin{Le}
  \label{le:1}
  For any $p>0$, there exists $C_p>0$ such that, with probability at
  least $1-|\Lambda|^{-p}$, if $E_j(\omega,\Lambda)$ is in the
  localized regime i.e. if $E_j(\omega,\Lambda)\in I$ then the
  diameter of $C(\varphi_j(\omega,\Lambda))$ is less than
  $C_p\log|\Lambda|$.
\end{Le}
\noindent Hence, in the localized regime, localization centers for an
eigenfunction can be at most as far as $C\log|\Lambda|$ from each
other. From now on, a localization center for a function $\varphi$
will denote any point in the set of localization centers $C(\varphi)$
and let $x_j(\omega,\Lambda)$ be a localization center for
$\varphi_j(\omega,\Lambda)$.\\
\subsection{Uniform Poisson convergence for the joint (energy,center)-distribution}
\label{sec:strong-poiss-conv-1}
We now place ourselves in the same setting as in
section~\ref{sec:strong-poiss-conv}. We prove
\begin{Th}
  \label{thr:11}
  Assume (W), (M) and (Loc) hold. Let $E_0$ be an energy in $I$ such
  that $\nu(E_0)>0$.\\
  Then, there exists $\delta>0$, such that, 
  \begin{itemize}
  \item for any sequences of
    intervals $I_1=I^\Lambda_1,\dots,I_p=I^\Lambda_p$ in
    $|\Lambda|^{1-\beta}\cdot[-1,1]$ satisfying~\eqref{eq:32},
  \item for any sequences of cubes
    $C_1=C^\Lambda_1,\dots,C_p=C^\Lambda_p$ in $[-1/2,1/2]^d$
  \end{itemize}
  one has, for any sequences of integers
  $k_1=k^\Lambda_1,\cdots,k_p=k^\Lambda_p\in\N^p$,
  \begin{multline*}
    \lim_{|\Lambda|\to+\infty}
    \left|\P\left(\left\{\omega;\
        \begin{split}
          &\#\left\{n;
            \begin{split}
              \xi_n(\omega,\Lambda)&\in I_1 \\ x_n/L&\in
              C_1 \end{split}
          \right\}=k_1\\&\vdots\hskip3cm\vdots\\
          &\#\left\{n;
            \begin{split}
              \xi_n(\omega,\Lambda)&\in I_p \\ x_n/L&\in C_p
            \end{split}
          \right\}=k_p
        \end{split}
      \right\}\right)-
    \prod_{j=1}^pe^{-|I_j||C_j|}\frac{(|I_j||C_j|)^{k_j}}{k_j!}\right|=0
  \end{multline*}
  where $x_n(\omega)=x_n(\omega,\Lambda_L)$ is the localization center
  associated to the eigenvalue $E_n(\omega,\Lambda_L)=E_0+L^d
  \xi_n(\omega,\Lambda)$.
\end{Th}
\noindent This result generalizes the results
of~\cite{MR2299191,MR2255080}.
\subsection{Covariant scaling joint (energy,center)-distribution}
\label{sec:covar-scal-joint}
Fix a sequence of scales $\ell=(\ell_{\Lambda})_{\Lambda}$ such that
\begin{equation}
  \label{eq:14}
  \frac{\ell_{\Lambda}}{\log|\Lambda|}\vers_{|\Lambda|\to+\infty}+\infty
  \quad\text{ and }\quad \ell_{\Lambda}\leq |\Lambda|^{1/d}.
\end{equation}
Pick $E_0\in I$ so that $\nu(E_0)>0$. Consider the point process
\begin{equation*}
  \Xi^2_\Lambda(\xi,x;E_0,\ell) = \sum_{j=1}^N
  \delta_{\nu(E_0)(E_j(\omega,\Lambda)-E_0)\ell_\Lambda^d}(\xi)
  \otimes\delta_{x_j(\omega)/\ell_\Lambda}(x). 
\end{equation*}
The process is valued in $\R\times\R^d$; actually, if
$c\,\ell_{\Lambda}\geq|\Lambda|^{1/d}$, it is valued in
$\R\times(-c,c)^d$. Define
\begin{equation*}
  c_\ell=\lim_{|\Lambda|\to+\infty}|\Lambda|^{1/d}\ell_\Lambda^{-1}
  \in [1,+\infty].
\end{equation*}
We prove
\begin{Th}[\cite{Ge-Kl:10}]
  \label{thr:5}
  The point process $\Xi^2_\Lambda(\xi,x;E_0,\ell)$ converges weakly
  to a Poisson process on $\R\times(-c_\ell,c_\ell)^d$ with intensity
  measure the Lebesgue measure.
\end{Th}
\noindent In the case $\ell_{\Lambda}=|\Lambda|^{1/d}$, the result of
Theorem~\ref{thr:5} was obtained in~\cite{MR2299191} (see
also~\cite{MR2355568,MR2255080}). In general, we see that, once the
energies and the localization centers are scaled covariantly, the
convergence to a Poisson process is true at any scale that is
essentially larger than the localization width. The scaling is very
natural; it is the one prescribed by the Heisenberg uncertainty
principle: the more precision we require in the energy variable, the
less we can afford in the space variable. In this respect, the
energies behave like a homogeneous symbol of degree $d$. This is quite
different from what one has in the case of the Laplace operator.
\subsection{Non-covariant scaling joint (energy,center)-distribution}
\label{sec:non-covar-scal}
One can also study what happens when the energies and localization
centers are not scaled covariantly. Consider two sequences of scales,
say $\ell=(\ell_{\Lambda})_{\Lambda}$ and
$\ell'=(\ell'_{\Lambda})_{\Lambda}$. Pick $E_0\in I$ so that
$\nu(E_0)>0$. Consider the point process
\begin{equation*}
  \Xi^2_\Lambda(\xi,x;E_0,\ell,\ell') = \sum_{j=1}^N
  \delta_{\nu(E_0)(E_j(\omega,\Lambda)-E_0)\ell_\Lambda^d}(\xi)
  \otimes\delta_{x_j(\omega)/\ell'_\Lambda}(x). 
\end{equation*}
Then, one proves
\begin{Th}[\cite{Ge-Kl:10}]
  \label{thr:7}
  Assume the sequences of increasing scales
  $\ell=(\ell_{\Lambda})_{\Lambda}$ and
  $\tilde\ell=(\tilde\ell_{\Lambda})_{\Lambda}$ satisfy~\eqref{eq:14}.
  Assume that
  \begin{equation}
    \label{eq:70}
    \text{if }\ell_L=o(L)\text{ then
    }\frac{\ell_{\Lambda_{L+\ell_L}}}{\ell_{\Lambda_L}}
    \vers_{|\Lambda|\to+\infty}1\text{ and }
    \frac{\ell'_{\Lambda_{L+\ell_L}}}{\ell'_{\Lambda_L}}
    \vers_{|\Lambda|\to+\infty}1.
  \end{equation}
  Let $J$ and $C$ be bounded measurable sets respectively in $\R$ and
  $(-c_{\tilde\ell},c_{\tilde\ell})^d\subset\R^d$. One has
  \begin{enumerate}
  \item if, for some $\rho>0$, one has
    $\displaystyle\frac{\tilde\ell_\Lambda}{\ell'_\Lambda}
    \leq|\Lambda|^{-\rho}$, then $\omega$-almost surely, for
    $\Lambda$ sufficiently large, 
    \begin{equation*}
      \int_{J\times C}\Xi^2_\Lambda(\xi,x;E_0,\ell,\tilde\ell)d\xi
      dx=0.
    \end{equation*}
  \item if, for some $\rho>0$, one has
  $\displaystyle\frac{\tilde\ell_\Lambda}{\ell'_\Lambda}
  \geq|\Lambda|^{\rho}$, then $\omega$-almost surely,
    \begin{equation*}
      \left(\frac{\ell_\Lambda}{\tilde\ell_\Lambda}\right)^{-d}
      \int_{J\times C}\Xi^2_\Lambda(\xi,x;E_0,\ell,\tilde\ell)d\xi
      dx\vers_{|\Lambda|\to+\infty}|J|\cdot|C|.
    \end{equation*}
  \end{enumerate} 
\end{Th}
\noindent Theorem~\ref{thr:7} proves that the local energy levels and
the localization centers become uniformly distributed in large energy
windows if one conditions the localization centers to a cube of much
smaller side-length. On the other hand, for a typical sample, if one
looks for eigenvalues in an energy interval much smaller than the
correctly scaled one with localization centers in a cube, then,
asymptotically, there are none.\\
Under assumption~\eqref{eq:14}, if one replaces the polynomial growth
or decay conditions on the ratio of scales by the condition that they
tend to $0$ or $\infty$, or if one omits condition~\eqref{eq:70}, the
results stays valid except for the fact that the convergence is not
almost sure anymore but simply holds in some $L^p$ norm.
\subsection{The level spacing statistics}
\label{sec:level-statistics}
\noindent Our goal is now to understand the level spacing statistics
for eigenvalues near $E_0\in I$. Pick $I_\Lambda$ a compact interval
containing $E_0$ such that its Lebesgue measure $|I_\Lambda|$ stays
bounded.\\
First, let us note that, by the existence of the density of states and
also Theorem~\ref{thr:3}, if $\nu(E_0)>0$, the spacing between
eigenvalues of $H_\omega(\Lambda)$ near $E_0$ is of size
$|\Lambda|^{-1}$. Hence, to study the statistics of level spacings in
$I_\Lambda$, $I_\Lambda$ should contain asymptotically infinitely many
energy levels of $H_{\omega,\Lambda}$. Let us study the number of
these levels.\\
\subsubsection{A large deviation principle for the eigenvalue counting
  function}
\label{sec:large-devi-princ}
Define the random numbers
\begin{equation}
  \label{eq:26}
  N(I_\Lambda,\omega,\Lambda):=\#\{j;\ E_j(\omega,\Lambda)\in
  I_\Lambda\}.
\end{equation}
Write $I_\Lambda=[a_\Lambda,b_\Lambda]$. We show that
$N(I_\Lambda,\omega,\Lambda)$ satisfies a large deviation principle
\begin{Th}[\cite{Ge-Kl:10}]
  \label{thr:16}
  Define $N(I_\Lambda)=N(b_\Lambda)-N(a_\Lambda)$ and assume that, for
  some $\nu\in(0,1)$,
  \begin{equation}
    \label{eq:1}
    N(I_\Lambda)\geq  |I_\Lambda|^{2-\nu}
  \end{equation}
  There exists $\delta>0$ and a sequence
  $(\varepsilon_\Lambda)_\Lambda$ such that, for
  $\varepsilon_\Lambda>0$, $\varepsilon_\Lambda\to0$ and one has
  \begin{equation}
    \label{eq:58}
    \P\left(\big|N(I_\Lambda,\omega,\Lambda)-N(I_\Lambda)|\Lambda|\big|\geq
      \varepsilon_\Lambda N(I_\Lambda)|\Lambda|\right)\leq
    e^{-(N(I_\Lambda)|\Lambda|)^\delta /\delta}.
  \end{equation}
\end{Th}
\noindent The large deviation principle~\eqref{eq:58} is meaningful
only if $N(I_\Lambda)|\Lambda|\to+\infty$; as $N$ is Lipschitz
continuous as a consequence of (W), this implies that
\begin{equation*}
  |\Lambda|\cdot|I_\Lambda|\vers
  +\infty\quad\text{ when }\quad|\Lambda|\to+\infty .
\end{equation*}
In this case, if $N(I_\Lambda)|\Lambda|$ satisfies~\eqref{eq:1}, one has
\begin{equation*}
  \esp(N(I_\Lambda,\omega,\Lambda))=N(I_\Lambda)|\Lambda|
  +o\left(N(I_\Lambda)|\Lambda|\right).
\end{equation*}
So~\eqref{eq:58} also says
\begin{equation*}
  \P\left(|N(I_\Lambda,\omega,\Lambda)-\esp(N(I_\Lambda,\omega,\Lambda))|\geq
    \varepsilon_\Lambda\, \esp(N(I_\Lambda,\omega,\Lambda))\right)\leq
  e^{-\esp(N(I_\Lambda,\omega,\Lambda))^\delta /\delta}.
\end{equation*}
\begin{Rem}
  \label{rem:1}
  Notice that the condition~\eqref{eq:1} allows for $I_\Lambda$
  to be centered at a point $E_0$ where $\nu(E_0)=0$ as long as the
  rate of vanishing of $\nu$ near $E_0$ is not too fast.\\
  Actually, all the results presented in this paper can be extended to
  this setting i.e. in
  Theorems~\ref{thr:3},~\ref{thr:2},~\ref{thr:4},~\ref{thr:18},~\ref{thr:8},~\ref{thr:11},~\ref{thr:5},~\ref{thr:7},~\ref{thr:1},~\ref{thr:DCS}
  and~\ref{thr:10}, one can replace the assumption $\nu(E_0)>0$
  by~\eqref{eq:1} (see\cite{Ge-Kl:10}). Of course, for the results to
  remain valid, in the definition of the points processes or the
  empirical distributions, one has to replace the normalization
  constant $|\Lambda|\nu(E_0)$ by $|\Lambda|\,
  N(I_\Lambda)/|I_\Lambda|$.
\end{Rem}
\subsubsection{The level spacing statistics near a given energy}
\label{sec:level-spac-stat}
Define $\mathcal{E}$ to be the set of energies $E$ such that
$\nu(E)=N'(E)$ exists and
\begin{equation*}
  \lim_{|x|+|y|\to0}\frac{N(E+x)-N(E+y)}{x-y}=\nu(E).
\end{equation*}
The requirement on the points in $\mathcal{E}$ is somewhat stronger
than asking for the simple existence of $\nu(E)$. Nevertheles, one
proves that the set $\mathcal{E}$ is of full Lebesgue measure.  It
clearly contains the continuity points of $\nu(E)$.
\vskip.1cm\noindent Fix $E_0\in\mathcal{E}$. If
$I_\Lambda=[a_\Lambda,b_\Lambda]$ is such that
$\displaystyle\sup_{I_\Lambda}|x|\vers_{|\Lambda|\to+\infty}0$, then
\begin{equation*}
  N_\Lambda(E_0+I_\Lambda)=\nu(E_0)|I_\Lambda||\Lambda|(1+o(1))
  \quad\text{as}\quad  |\Lambda|\to+\infty. 
\end{equation*}
Consider the renormalized eigenvalue spacings: for $1\leq j\leq N$,
\begin{equation*}
  \delta E_j(\omega,\Lambda)=|\Lambda|\,\nu(E_0)(E_{j+1}(\omega,\Lambda)-
  E_j(\omega,\Lambda))\geq0.
\end{equation*}
Define the empirical distribution of these spacings to be the random
numbers, for $x\geq0$
\begin{equation*}
  DLS(x;I_\Lambda,\omega,\Lambda)=\frac{\#\{j;\ E_j(\omega,\Lambda)\in
    I_\Lambda,\ \delta E_j(\omega,\Lambda)\geq x\}
  }{N(I_\Lambda,\omega,\Lambda)}.
\end{equation*}
We first study the level spacings distributions of the energies
inside an interval that shrink to a point.\\
We prove
\begin{Th}[\cite{Ge-Kl:10}]
  \label{thr:1}
  Fix $E_0\in\mathcal{E}$ such that $\nu(E_0)>0$ and pick
  $(I_\Lambda)_\Lambda$ a sequence of intervals centered at $E_0$ such
  that
  $\displaystyle\sup_{I_\Lambda}|x|\vers_{|\Lambda|\to+\infty}0$.\\
  Assume that, for some $\delta>0$, one has
  \begin{equation}
    \label{eq:43}
    |\Lambda|^{1-\delta}\cdot|I_\Lambda|
    \vers_{|\Lambda|\to+\infty}    +\infty
    \quad \text{ and }\quad \text{if }\ell_L=o(L)\text{ then
    }\frac{|I_{\Lambda_{L+\ell_L}}|}{|I_{\Lambda_L}|}
    \vers_{L\to+\infty}1.
  \end{equation}
  Then, with probability $1$ , as $|\Lambda|\to+\infty$,
  $DLS(x;I_\Lambda,\omega,\Lambda)$ converges uniformly to the
  distribution $x\mapsto e^{-x}$, that is, with probability $1$,
  \begin{equation*}
    \sup_{x\geq0}\left|DLS(x;I_\Lambda,\omega,\Lambda)
          -e^{-x}\right|\vers_{|\Lambda|\to+\infty}0.
  \end{equation*}
\end{Th}
\noindent Hence, the rescaled level spacings behave as if the
eigenvalues were i.i.d. uniformly distributed random variables (see
\cite{MR0070874} or section 7 of~\cite{MR0216622}). This distribution
for the level spacings is the one predicted by physical heuristics in
the localized regime
(\cite{citeulike:693492,RevModPhys.57.287,citeulike:3832118,Th:74}). It
is also in accordance with Theorem~\ref{thr:3}.
In~\cite{MR97d:82046,CGK:10}, the domains in energy where the
statistics could be studied were much smaller than the ones considered
in Theorem~\ref{thr:1}. Indeed, the energy interval was of order
$|\Lambda|^{-1}$ whereas, here, it is assumed to tend to $0$ but be
large when compared to $|\Lambda|^{-1}$. In particular,
in~\cite{MR97d:82046,CGK:10}, the intervals were not large enough to
enable the computation of statistics of levels as not enough levels
were involved: the intervals typically contained only finitely many
intervals.
\par The first condition in~\eqref{eq:43} ensures that $I_\Lambda$
contains sufficiently many eigenvalues of $H_\omega(\Lambda)$. The
second condition in~\eqref{eq:43} is a regularity condition of the
decay of $|I_\Lambda|$. If one omits either or both of these two
conditions and only assumes that
$|\Lambda|\cdot|I_\Lambda|\to+\infty$, one still gets convergence in
probability of $DLS(x;I_\Lambda,\omega,\Lambda)$ to $e^{-x}$ i.e.
\begin{equation*}
  \P\left(\sup_{x\geq0}\left|DLS(x;I_\Lambda,\omega,\Lambda)
      -e^{-x}\right|\geq\varepsilon\right)\vers_{|\Lambda|\to+\infty}0.
\end{equation*}
\subsubsection{The level spacing statistics on macroscopic energy intervals}
\label{sec:level-spac-stat-1}
Theorem~\ref{thr:1} seems optimal as the density of states at $E_0$
enters into the correct rescaling to obtain a universal result. Hence,
the distribution of level spacings on larger intervals needs to take
into account the variations of the density of states on these
intervals. Indeed, on intervals of non vanishing size, we compute the
asymptotic distribution of the level spacings when one omits the local
density of states in the spacing and obtain
\begin{Th}[\cite{Ge-Kl:10}]
  \label{thr:9}
  Pick $J\subset I$ a compact interval such
  $\lambda\mapsto\nu(\lambda)$ be continuous on $J$ and
  $N(J):=\int_J\nu(\lambda)d\lambda>0$. Define the renormalized
  eigenvalue spacings, for $1\leq j\leq N$,
  \begin{equation*}
    \delta_J E_j(\omega,\Lambda)=|\Lambda|N(J)
    (E_{j+1}(\omega,\Lambda)-E_j(\omega,\Lambda))\geq0
  \end{equation*}
  and the empirical distribution of these spacing to be the random
  numbers, for $x\geq0$
  \begin{equation*}
    DLS'(x;J,\omega,\Lambda)=\frac{\#\{j;\ E_j(\omega,\Lambda)\in
      J,\ \delta_J E_j(\omega,\Lambda)\geq x\}
    }{N(J,\omega,\Lambda)}.
  \end{equation*}
  Then, as $|\Lambda|\to+\infty$, with probability 1,
  $DLS'(x;J,\omega,\Lambda)$ converges uniformly to the distribution
  $x\mapsto g_{\nu,J}(x)$ where
  \begin{equation}
    \label{eq:31}
    g_{\nu,J}(x)=\int_{J}\e^{-\nu_J(\lambda)x}\nu_J(\lambda)
    d\lambda\text{ where }\nu_J=\frac{1}{N(J)}\nu.
  \end{equation}
\end{Th}
\noindent We see that, in the large volume limit, the rescaled
level spacings behave as if the eigenvalues were i.i.d. random
variables distributed according to the density
$\frac{1}{N(J)}\nu(\lambda)$ i.e. to the density of states
normalized to be a probability measure on $J$ (see section 7
of~\cite{MR0216622}).\\
In Theorem~\ref{thr:9}, we assumed the density of states to be
continuous. This is known to hold in the large coupling limit if the
density of the distribution of the random variables is sufficiently
smooth (see~\cite{MR89c:82050}).
\section{The localization center spacing statistics}
\label{sec:center-statistics}
Pick $E_0$ as above. Inside the cube $\Lambda$, the number of centers
that corresponds to energies in $I_\Lambda$ is roughly equal to
$\nu(E_0) |I_\Lambda| N$. Thus, if we assume that the localization
centers are uniformly distributed as is suggested by
Theorems~\ref{thr:5} and~\ref{thr:7}, the reference mean spacing
between localization centers is of size $(|\Lambda|/(\nu(E_0)
|I_\Lambda| |\Lambda|)^{1/d} = (\nu(E_0)|I_\Lambda|)^{-1/d}$. This
motivates the following definition.\\
Define the empirical distribution of center spacing to be the random
number
\begin{equation}
  \label{DCS}
  DCS(s;I_\Lambda,\omega,\Lambda)=\frac{\#\left\{j;\ E_(\omega,\Lambda)\in
      I_\Lambda,\displaystyle\  \frac{\min_{i\neq j}|x_j(\omega)-x_i(\omega)|}
      {\sqrt[d]{\nu(E_0)|I_\Lambda|}}\geq s \right\} }
  {N(I_\Lambda,\omega,\Lambda)}
\end{equation}
where $N(I_\Lambda,\omega,\Lambda)$ is defined in~\eqref{eq:26}.\\
We prove an analogue of Theorems~\ref{thr:1}, namely
\begin{Th}[\cite{Ge-Kl:10}]
  \label{thr:DCS}
  Pick $E_0\in I$ such that $\nu(E_0)>0$.  Assume (H.3) and
  \begin{equation*}
    |I_\Lambda|=o\left(\frac1{\log^d|\Lambda|}\right).
  \end{equation*}
  Then, as $|\Lambda|\to+\infty$, in probability,
  $DCS(s;I_\Lambda,\omega,\Lambda)$ converges uniformly to the
  distribution $x\mapsto e^{-s^d}$, that is, for any $\varepsilon>0$,
  \begin{equation*}
    \P\left(\left\{\omega;\
        \sup_{s\geq0}\left|DCS(s;I_\Lambda,\omega,\Lambda)-e^{-s^d} 
        \right|\geq\varepsilon\right\}\right)\vers_{\Lambda\nearrow\R^d} 0.
  \end{equation*}
\end{Th}
\noindent Of course, Theorem~\ref{thr:9} also has an analogue for
localization centers.
\section{Another point of view}
\label{sec:another-point-view}
In the present section, we want to adopt a different point of view on
the spectral statistics. Instead of discussing the statistics of the
eigenvalues of the random system restricted to some finite box in the
large box limit, we will describe the spectral statistics of the
infinite system in the localized phase. Let $I$ be an interval in the
region of complete localization. Then, it is well known
(\cite{MR2509110,MR2002h:82051,MR2203782}) that, in this region, the
following property holds
\begin{description}
\item[(Loc')] there exists $\gamma>0$ such that, with probability $1$,
  if $E\in I\cap\sigma(H_\omega)$ and $\varphi$ is a normalized
  eigenfunction associated to $E$ then, for $x(E)\in\Z^d$, a maximum
  of $x\mapsto\|\varphi\|_x$, for some $C_\omega>0$, one has, for
  $x\in\R^d$,
    \begin{equation*}
      |\varphi(x)|\leq C_\omega (1+|x(E)|^2)^{q/2}e^{-\gamma
        |x-x(E)|};
    \end{equation*}
    moreover, one has $\esp(C_\omega)<+\infty$.
\end{description}
As above $x(E)$ is called {\it a center of localization} for the
energy $E$ or for the associated  eigenfunction $\varphi$.\\
Without restriction on generality, we assume that
$\sigma(H_\omega)\cap I=I$ $\omega$-almost surely. Hence, any
sub-interval of $I$ contains infinitely many eigenvalues and to define
statistics, we need a way to enumerate these eigenvalues. To do this,
we use the localization centers; namely, we prove
\begin{Pro}[\cite{Ge-Kl:10}]
  \label{pro:1}
  Fix $q>2d$. Then, there exists $\gamma>0$ such that, $\omega$-almost
  surely, there exists $C_\omega>1$ such that
  \begin{enumerate}
  \item if $x(E)$ and $x'(E)$ are two centers of localization for
    $E\in I$ then
    \begin{equation*}
      |x(E)-x'(E)|\leq \gamma^{-2}(\log\langle x(E)\rangle+\log
      C_\omega)^{1/\xi}. 
    \end{equation*}
  \item for $L\geq1$, pick $I_L\subset I$ such that
    $L^dN(I_L)\to+\infty$ (see Theorem~\ref{thr:16}); if $N(I_L,L)$
    denotes the number of eigenvalues of $H_\omega$ having a center of
    localization in $\Lambda_L$, then
    \begin{equation*}
      N(I_L,L)=N(I_L)\,|\Lambda_L|\,(1+o(1)).
    \end{equation*}
  \end{enumerate}
\end{Pro}
\noindent Point (1) is proved in~\cite{MR2203782} (see Corollary 3 and
its proof). Points (2) and (3) are proved in~\cite{Ge-Kl:10}.\\
For $L\geq1$, pick $I_L\subset I$ such that
$L^d\,N(I_L)\to+\infty$. In view of Proposition~\ref{pro:1}, we can
consider the level spacings for the eigenvalues of $H_\omega$ having a
localization center in $\Lambda_L$; indeed, for $L$ large, there are
only finitely many such eigenvalues, let us enumerate them as
$E_1(\omega,L)\leq E_2(\omega,L)\leq \cdots\leq E_N(\omega,L)$ where
we repeat them according to multiplicity. Consider the renormalized
eigenvalue spacings, for $1\leq j\leq N$,
\begin{equation*}
  \delta E_j(\omega,L)=|\Lambda_L|\,(E_{j+1}(\omega,L)-
  E_j(\omega,L))\geq0.
\end{equation*}
Define the empirical distribution of these spacing to be the random
numbers, for $x\geq0$
\begin{equation*}
  DLS(x;I_L,\omega,L)=\frac{\#\{j;\ E_j(\omega,L)\in
    I_L,\ \delta E_j(\omega,L)\geq x\}
  }{N(I_L,L)}
\end{equation*}
Then, we prove
\begin{Th}[\cite{Ge-Kl:10}]
  \label{thr:10}
  One has
  \begin{itemize}
  \item if $E_0\in\mathcal{E}\cap I_L$ s.t. $\nu(E_0)>0$ and
    $|I_L|\to0$ and satisfies~\eqref{eq:43}, then, $\omega$-almost
    surely, for $x\geq0$
    \begin{equation*}
      \lim_{L\to+\infty}\sup_{x\geq0}
      \left|DLS(x;I_L,\omega,L)-e^{-\nu(E_0)x}\right|=0;
    \end{equation*}
  \item if, for all $L$ large, $I_L=J$ such that $\nu(J)>0$ and $\nu$
    is continuous on $J$ then, $\omega$-almost surely, one has
    \begin{equation*}
      \lim_{L\to+\infty}\sup_{x\geq0}\left|DLS(x;I_L,\omega,L)
        -g_{\nu,J}(N(J)\,x)\right|=0
    \end{equation*}
    where $g_{\nu,J}$ is defined in~\eqref{eq:31}.
  \end{itemize}
\end{Th}
\noindent In the first part of Theorem~\ref{thr:10}, if~\eqref{eq:43}
is not satisfied, then the convergence still holds in probability.

%

\begin{thebibliography}{10}

\bibitem{MR2002h:82051}
Michael Aizenman, Jeffrey~H. Schenker, Roland~M. Friedrich, and Dirk
  Hundertmark.
\newblock Finite-volume fractional-moment criteria for {A}nderson localization.
\newblock {\em Comm. Math. Phys.}, 224(1):219--253, 2001.
\newblock Dedicated to Joel L. Lebowitz.

\bibitem{MR89c:82050}
A.~Bovier, M.~Campanino, A.~Klein, and J.~F. Perez.
\newblock Smoothness of the density of states in the {A}nderson model at high
  disorder.
\newblock {\em Comm. Math. Phys.}, 114(3):439--461, 1988.

\bibitem{CGK:10}
Jean-Michel Combes, Fran{{\c c}}ois Germinet, and Abel Klein.
\newblock Poisson statistics for eigenvalues of continuum random
  {Schr{\"o}dinger} operators.
\newblock Preprint.

\bibitem{MR2505733}
Jean-Michel Combes, Fran{{\c c}}ois Germinet, and Abel Klein.
\newblock Generalized eigenvalue-counting estimates for the {A}nderson model.
\newblock {\em J. Stat. Phys.}, 135(2):201--216, 2009.

\bibitem{MR2203782}
Francois Germinet and Abel Klein.
\newblock New characterizations of the region of complete localization for
  random {S}chr{\"o}dinger operators.
\newblock {\em J. Stat. Phys.}, 122(1):73--94, 2006.

\bibitem{Ge-Kl:10}
Fran{{\c c}}ois Germinet and Fr{\' e}d{\'e}ric Klopp.
\newblock Spectral statistics for random operators in the localized regime.
\newblock in progress.

\bibitem{citeulike:693492}
Martin Janssen.
\newblock Statistics and scaling in disordered mesoscopic electron systems.
\newblock {\em Physics Reports}, 295(1-2):1--91, March 1998.

\bibitem{MR2299191}
Rowan Killip and Fumihiko Nakano.
\newblock Eigenfunction statistics in the localized {A}nderson model.
\newblock {\em Ann. Henri Poincar{\'e}}, 8(1):27--36, 2007.

\bibitem{MR2509110}
Werner Kirsch.
\newblock An invitation to random {S}chr{\"o}dinger operators.
\newblock In {\em Random {S}chr{\"o}dinger operators}, volume~25 of {\em Panor.
  Synth{\`e}ses}, pages 1--119. Soc. Math. France, Paris, 2008.
\newblock With an appendix by Fr{\'e}d{\'e}ric Klopp.

\bibitem{Kl:10}
Fr{\'e}d{\'e}ric Klopp.
\newblock Decorrelation estimates for the discrete anderson model.
\newblock Preprint arXiv: {\tt http://fr.arxiv.org/abs/1004.1261}, 2010.

\bibitem{RevModPhys.57.287}
Patrick~A. Lee and T.~V. Ramakrishnan.
\newblock Disordered electronic systems.
\newblock {\em Rev. Mod. Phys.}, 57(2):287--337, Apr 1985.

\bibitem{MR97d:82046}
Nariyuki Minami.
\newblock Local fluctuation of the spectrum of a multidimensional {A}nderson
  tight binding model.
\newblock {\em Comm. Math. Phys.}, 177(3):709--725, 1996.

\bibitem{citeulike:3832118}
A.~Mirlin.
\newblock Statistics of energy levels and eigenfunctions in disordered systems.
\newblock {\em Physics Reports}, 326(5-6):259--382, March 2000.

\bibitem{MR2255080}
Fumihiko Nakano.
\newblock The repulsion between localization centers in the {A}nderson model.
\newblock {\em J. Stat. Phys.}, 123(4):803--810, 2006.

\bibitem{MR2355568}
Fumihiko Nakano.
\newblock Distribution of localization centers in some discrete random systems.
\newblock {\em Rev. Math. Phys.}, 19(9):941--965, 2007.

\bibitem{MR0216622}
Ronald Pyke.
\newblock Spacings. ({W}ith discussion.).
\newblock {\em J. Roy. Statist. Soc. Ser. B}, 27:395--449, 1965.

\bibitem{Th:74}
David Thouless.
\newblock Electrons in disordered systems and the theory of localization.
\newblock {\em Physical Reports}, 13:93--142, 1974.

\bibitem{MR0070874}
Lionel Weiss.
\newblock The stochastic convergence of a function of sample successive
  differences.
\newblock {\em Ann. Math. Statist.}, 26:532--536, 1955.

\end{thebibliography}
%
\def\cprime{$'$} \def\cydot{\leavevmode\raise.4ex\hbox{.}} \def\cprime{$'$}

\end{document}